\documentclass[12pt,a4paper]{article}
\title{A Bound for the Number of Different Basic Solutions 
Generated by the Simplex Method}
\author{Tomonari Kitahara\footnote{Graduate School of Decision Science and Technology, Tokyo Institute of Technology, 
2-12-1-W9-62, Ookayama, Meguro-ku, Tokyo, 152-8552, Japan. 
Tel.: +81-3-5734-2896
Fax: +81-3-5734-2947 
E-mail: kitahara.t.ab@m.titech.ac.jp}\hspace{1ex}
 and Shinji Mizuno\footnote{Graduate School of Decision Science and Technology, Tokyo Institute of Technology, 
2-12-1-W9-58, Ookayama, Meguro-ku, Tokyo, 152-8552, Japan. 
Tel.: +81-3-5734-2816
Fax: +81-3-5734-2947 
E-mail: mizuno.s.ab@m.titech.ac.jp}}
\date{October 15, 2010; minor revision December 1, 2010}
\newtheorem{lemma}{Lemma}[section]
\newtheorem{theorem}{Theorem}[section]
\newtheorem{corollary}{Corollary}[section]
\usepackage{amsmath,amssymb}
\newcommand{\qed}{\hfill$\blacksquare$\par}
\begin{document}
\maketitle
\begin{abstract}
In this short paper, we give an upper bound for the number of different basic feasible solutions 
generated by the simplex method for linear programming problems having optimal solutions. 
The bound is  polynomial of the number of constraints, the number of variables, and 
the ratio between  the minimum and the maximum values of all the positive elements of primal basic feasible solutions.  
When the primal problem is nondegenerate, it becomes a bound for the number of iterations. 
We show some basic results when it is applied to special linear programming problems. 
The results include strongly polynomiality of the simplex method for Markov Decision Problem by Ye \cite{ye} 
and utilize its analysis.
\end{abstract}
\section{Introduction}
The simplex method for solving linear programming problems (LP)
was originally developed by Dantzig \cite{dantzig}. 
The simplex method works very efficiently in practice and it has been widely used 
for years. 
In spite of the practical efficiency of the simplex method, 
we do not have any good bound for the number of iterations.  
Klee and Minty \cite{km} showed that the simplex method needs an exponential number of iterations 
for an elaborately designed problem. 

We analyze the primal simplex method with the most negative pivoting rule (Dantzig's rule) or 
the best improvement pivoting rule under the condition that the primal problem has an optimal solution. 
We give an upper bound for the number of different 
basic feasible solutions (BFSs) generated by the simplex method.
The bound is 
\[
n \lceil m\frac{\gamma}{\delta}\log (m\frac{\gamma}{\delta})\rceil,
\]
where $m$ is the number of constraints, $n$ is the number of variables, 
$\delta$ and $\gamma$ are the minimum and the maximum 
values of all the positive elements of primal BFSs, respectively, 
and $\lceil a \rceil$ is the smallest integer bigger than $a\in\Re$.   
When the primal problem is nondegenerate, it becomes a bound for the number of iterations. 
Note that the bound depends only on the constraints of LP, but not the objective function.

Our work is motivated by a recent research by Ye \cite{ye}. 
He shows that the simplex method is strongly polynomial for the
Markov Decision Problem. 
We apply the analysis in \cite{ye} to general LPs and obtain the upper bound. 
Our results include his strong polynomiality.  

When we apply our result to an LP where a constraint matrix is totally unimodular 
and a constant vector $b$ of constraints is integral, the number of different solutions generated by the simplex method 
is at most
\[n\lceil m\|b\|_1\log (m\|b\|_1)\rceil.
\]

The paper is organized as follows. 
In section 2, we explain basic notions of an LP and 
we briefly review the simplex method. 
In section 3, analyses of the simplex method are conducted to show our results. 
In section 4, applications of our results to special LPs are discussed.

\section{The Simplex Method for LP}
In this paper, we consider the linear programming problem of the standard form 
\begin{equation}
\begin{array}{ll}
\min & c^Tx,\\
{\rm subject~ to}& Ax=b,~ x\ge 0, \label{primal}
\end{array}
\end{equation}
where $A\in\Re^{m\times n},~b\in\Re^m$ and $c\in\Re^n$ are given data,
and $x\in\Re^n$ is a variable vector. 
The dual problem of (\ref{primal}) is  
\begin{equation}
\begin{array}{ll}
\max&b^Ty,\\
{\rm subject~to}&A^Ty+s=c,~s\ge 0, \label{dual}
\end{array}
\end{equation}
where $y\in\Re^m$ and $s\in\Re^n$ are variable vectors.\\

We assume that ${\rm rank}(A)=m$, the primal problem (\ref{primal}) has an optimal solution 
and an initial BFS $x^0$ is available. 
Let $x^*$ be an optimal basic feasible solution of (\ref{primal}), 
$(y^*,s^*)$ be an optimal solution of (\ref{dual}), 
and $z^*$ be the optimal value of (\ref{primal}) and (\ref{dual}). \\

Given a set of indices $B\subset \{1,2,\dots,n\}$, we split the constraint matrix $A$, 
the objective vector $c$, and the variable vector $x$ according to $B$ and $N=\{1,2,\dots,n\}-B$ like
\[
A=[A_B,A_N],~c=
\left[
\begin{array}{c}
c_B\\
c_N
\end{array}
\right],~
x=
\left[
\begin{array}{c}
x_B\\
x_N
\end{array}
\right].
\]
Define the set of bases 
\[
{\cal B}=\{B\subset\{1,2,\dots,n\}|~|B|=m,~\det(A_B)\not=0\}.
\]
Then a primal basic feasible solution for $B\in {\cal B}$ and $N=\{1,2,\dots,n\}-B$ is written as
\[
x_B=A_B^{-1}b\ge 0,~x_N=0.
\]
Let $\delta$ and $\gamma$ be the minimum and the maximum values of all the positive elements of BFSs. 
Hence for any BFS $\hat{x}$ and any $j\in\{1,2,\dots,n\}$, if $\hat{x}_j\not=0$, we have 
\begin{equation}
\delta\le \hat{x}_j \le \gamma. \label{l-u}
\end{equation}
Note that these values depend only on $A$ and $b$, but not on $c$.\\ 

Let $B^t\in{\cal B}$ be the basis of the $t$-th iteration of the simplex method
and set $N^t=\{1,2,\dots,n\}-B^t$. 
Problem (\ref{primal}) can be written as 
\[
\begin{array}{ll}
\min& c_{B^t}^TA_{B^t}^{-1}b+(c_{N^t}-A_{N^t}^T(A_{B^t}^{-1})^Tc_{B^t})^Tx_{N^t}, \\
{\rm subject~to}&x_{B^t}=A_{B^t}^{-1}b-A_{B^t}^{-1}A_{N^t}x_{N^t},\\
                &x_{B^t}\ge0,~x_{N^t}\ge0.
                
\end{array}
\]
The coefficient vector $\bar{c}_{N^t}=c_{N^t}-A_{N^t}^T(A_{B^t}^{-1})^Tc_{B^t}$ 
is called a reduced cost vector. 
When $\bar{c}_{N^t}\ge 0$, the current solution is optimal. 
Otherwise we conduct a pivot. 
That is, we choose one nonbasic variable (entering variable) and increase the variable 
until one basic variable (leaving variable) becomes zero.
Then we exchange the two variables. 
Several rules for choosing the entering variable have been proposed. 
For example, the most negative rule, the best improvement rule, and the minimum index rule 
are well known. 
Under the most negative rule, we choose a nonbasic variable whose reduced cost is minimum.
To put it precisely, 
we choose an index 
\[j^t_{MN}=\arg\min_{j\in N_t}\bar{c}_j.\]
We set $\Delta^t=-\bar{c}_{j^t_{MN}}$. 

In the case of the best improvement rule, we choose an entering variable so that 
the objective value at the next iterate is minimum. 
We summarize the notations in Table \ref{notations}. 
\begin{table}
\caption{Notations}
\label{notations}
\begin{center}
\begin{tabular}{|cll|}
\hline
$x^*$&:&an optimal basic feasible solution of (\ref{primal})\\
$(y^*,s^*)$&:&an optimal solution of (\ref{dual}) \\
$z^*$&:&the optimal value of (\ref{primal})\\
$x^t$&:&the $t$-th iterate of the simplex method\\
$B^t$&:&the basis of $x^t$\\ 
$N^t$&:&the nonbasis of $x^t$\\
$\bar{c}_{N^t}$&:&the reduced cost vector at $t$-th iteration\\
$\Delta^t$&:&$-\min_{j\in N_t}\bar{c}_j$\\
$j^t_{MN}$&:& an index chosen by the most negative rule at $t$-th iteration\\
\hline
\end{tabular}
\end{center}
\end{table}

\section{Analysis of the Simplex Method}
Our analysis is motivated by a recent work by Ye \cite{ye}, 
where he shows a strongly polynomial result of the simplex method for the Markov Decision Problem, a special class of LP. 
We apply his analysis to general LPs and obtain an upper bound for the number of 
different basic feasible solutions generated by the simplex method. 
Later we confirm that our results include his strongly polynomiality. \\

In the next lemma, we get a lower bound of the optimal value at each iteration of the simplex method. 
\begin{lemma}
Let $z^*$ be the optimal value of Problem (\ref{primal}) and $x^t$ be the $t$-th iterate generated by the simplex method 
with the most negative rule. 
Then we have 
\begin{equation}
z^* \ge c^Tx^t-\Delta^t m\gamma. \label{ie1}
\end{equation}
\end{lemma}
{\bf Proof.} Let $x^*$ be a basic optimal solution of Problem (\ref{primal}). 
Then we have 
\[
\begin{array}{ll}
z^* &= c^Tx^*\\
    &= c_{B^t}^TA_{B^t}^{-1}b+\bar{c}_{N^t}^Tx^*_{N^t}\\
    &\ge c^Tx^t-\Delta^t e^Tx^*_{N^t}\\
    &\ge c^Tx^t-\Delta^t m \gamma,
\end{array}
\]
where the second inequality follows since $x^*$ has at most $m$ positive elements and 
each element is bounded above by $\gamma$. 
Thus we get the inequality (\ref{ie1}).\qed
\vspace{5mm}

Next we show a constant reduction of the gap between the objective function value 
and the optimal value at each iteration when 
an iterate is updated. 
The result is interesting becasuse the reduction rate $(1-\frac{\delta}{m\gamma})$ does not depend on the objective vector $c$.
\begin{theorem}\label{reduction}
Let $x^t$ and $x^{t+1}$ be the $t$-th and $(t+1)$-th iterates 
generated by the simplex method with the most negative rule.   
If $x^{t+1}\not = x^t$, then we have 
\begin{equation}
c^Tx^{t+1}-z^*\le (1-\frac{\delta}{m\gamma})(c^Tx^t-z^*) \label{ie2}.
\end{equation}
\end{theorem}
{\bf Proof.} Let $x^t_{j^t_{MN}}$ be the entering variable chosen at the $t$-th iteration.
If $x^{t+1}_{j^t_{MN}}=0$, then we have $x^{t+1}=x^t$, a contradiction occurs. 
Thus $x^{t+1}_{j^t_{MN}}\not=0$, and we have $x^{t+1}_{j^t_{MN}}\ge \delta$ 
from (\ref{l-u}).
Then we have 
\[
\begin{array}{ll}
c^Tx^t-c^Tx^{t+1} &= \Delta^t x^{t+1}_{j^t_{MN}}\\
                  &\ge \Delta^t \delta\\
                  &\ge \frac{\delta}{m\gamma}(c^Tx^t-z^*),
\end{array}
\]
where the last inequality comes from (\ref{ie1}).
The desired inequality readily follows from the above inequality.\qed
\vspace{5mm}

Note that under the best improvement pivoting rule, 
the objective function reduces at least as much as that with the most negative rule.   
Thus the next corollary easily follows. 
\begin{corollary} \label{best-improvement}
Let $x^t$ and $x^{t+1}$ be the $t$-th and $(t+1)$-th iterates 
generated by the simplex method with the best improvement rule.   
If $x^{t+1}\not = x^t$, then we also have (\ref{ie2}).
\end{corollary}
\vspace{5mm}

From Theorem \ref{reduction} and Corollary \ref{best-improvement},
we can easily get an upper bound for the number of different BFSs generated by 
the simplex method. 
\begin{corollary} \label{bound-using-object-function}
Let $\bar{x}$ be a second optimal BFS of LP (\ref{primal}), that is, 
a minimal BFS except for optimal BFSs. 
When we apply the simplex method with the most negative rule or 
the best improvement rule for LP (\ref{primal}) 
from an initial BFS $x^0$, we encounter at most 
\[
  \lceil m \frac{\gamma}{\delta} \log(\frac{c^T x^0 - z^*}{c^T \bar{x} - z^*}) \rceil 
\]
different BFSs.
\end{corollary}
{\bf Proof.} 
Let $x^t$ be the $t$-th iterates generated by the simplex method and 
let $\tilde{t}$ be the number of different BFSs appearing up to this iterate. 
Then we have 
\[
   c^T x^t - z^* \le ( 1- \frac{\delta}{m \gamma})^{\tilde{t}} ( c^T x^0 - z^*) 
\]
from (\ref{ie2}). If $\tilde{t}$ is bigger than or equal to the number in the corollary, 
we easily get  
\[
   c^T x^t - z^* < c^T \bar{x} - z^*.
\]
Since $\bar{x}$ is a second optimal BFS of LP (\ref{primal}), 
$x^t$ must be an optimal BFS from the inequality above. 
\qed
Note that the bound in the corollary above depends on the objective function. 
In the succeeding discussion, we will have another bound which is independent of 
the objective function.

The next Lemma states that if the current solution is not optimal, there is a basic variable 
which has an upper bound proportional to the gap between the objective value and the optimal value.  
\begin{lemma}\label{slack}
Let $x^t$ be the $t$-th iterate generated by the simplex method.
If $x^t$ is not optimal, there exists $\bar{j}\in B^t$ such that $x^t_{\bar{j}}>0$ and  
\[
s^*_{\bar{j}}\ge \frac{1}{m x^t_{\bar{j}}}(c^Tx^t-z^*), 
\]
where $s^*$ is an optimal slack vector of (\ref{dual}). 
Then for any $k$, the $k$-th iterate $x^k$ satisfies 
\[
x^k_{\bar{j}}\le \frac{m(c^Tx^k-z^*)}{c^Tx^t-z^*}x_{\bar{j}}^t.
\] 
\end{lemma}
{\bf Proof.} Since 
\[
c^Tx^t-z^*=(x^t)^Ts^*=\sum_{j\in B^t}x^t_js^*_j,
\]
there exists $\bar{j}\in B^t$ which satisfies 
\[
s^*_{\bar{j}}x^t_{\bar{j}}\ge \frac{1}{m}(c^Tx^t-z^*), 
\] 
or equivalently, $x^t_{\bar{j}}>0$ and
\[
s^*_{\bar{j}}\ge \frac{1}{mx^t_{\bar{j}}}(c^Tx^t-z^*).
\]
Moreover, for any $k$, the $k$-th iterate $x^k$ satisfies 
\[
c^Tx^k-z^*=(x^k)^Ts^*=\sum_{j=1}^nx^k_js^*_j\ge x^k_{\bar{j}}s^*_{\bar{j}},
\]
which implies 
\[
x^k_{\bar{j}}\le \frac{c^Tx^k-z^*}{s^*_{\bar{j}}}\le \frac{m(c^Tx^k-z^*)}{c^Tx^t-z^*}x^t_{\bar{j}}. 
\]
\qed

\begin{lemma}\label{index}
Let $x^t$ be the $t$-th iterate 
generated by the simplex method with 
the most negative rule or the best improvement rule. 
Assume that $x^t$ is not an optimal solution. 
Then there exists $\bar{j}\in B^t$ satisfying the following two conditions.
\begin{enumerate}
\item $x^t_{\bar{j}}>0$.
\item If the simplex method generates $\lceil m\frac{\gamma}{\delta}\log(m\frac{\gamma}{\delta})\rceil$ different basic feasible solutions 
after $t$-th iterate, 
then $x_{\bar{j}}$ becomes zero and stays zero.
\end{enumerate}
\end{lemma}
{\bf Proof.} For $k\ge t+1$, let $\tilde{k}$ be the 
number of different basic feasible solutions appearing 
between the $t$-th and $k$-th iterations. 
Then from Theorem \ref{reduction} and Lemma \ref{slack}, there exists $\bar{j}\in B_t$ which satisfies 
\[
x^k_{\bar{j}}\le m(1-\frac{\delta}{m\gamma})^{\tilde{k}}x^t_{\bar{j}}\le m\gamma(1-\frac{\delta}{m\gamma})^{\tilde{k}}.
\]
The second inequality follows from (\ref{l-u}).
Therefore, if $\tilde{k}> m\frac{\gamma}{\delta}\log(m\frac{\gamma}{\delta})$, 
we have $x^k_{\bar{j}}<\delta$, which implies $x^k_{\bar{j}}=0$ from the definition of 
$\delta$. \qed
\vspace{5mm}

The event described in Lemma \ref{index} can occur at most once for each variable. 
Thus we get the following result. 
\begin{theorem}
When we apply the simplex method with the most negative rule or 
the best improvement rule for LP (\ref{primal}) having optimal solutions, we encounter at most 
$n\lceil m\frac{\gamma}{\delta}\log(m\frac{\gamma}{\delta})\rceil$
different basic feasible solutions.
\end{theorem}
Note that the result is valid even if the simplex method fails to find an optimal solution because of a cycling.\\

If the primal problem is nondegenerate, we have $x^{t+1}\not= x^t$ for all $t$. 
This observation leads to a bound for the number of iterations of the simplex method.  
\begin{corollary}\label{iteration} 
If the primal problem is nondegenerate, the simplex method 
finds an optimal solution in at most $n\lceil m\frac{\gamma}{\delta}\log(m\frac{\gamma}{\delta})\rceil$ iterations.
\end{corollary}
\section{Applications to Special LPs}
\subsection{LP with a Totally Unimodular Matrix}
In this subsection, we consider an LP (\ref{primal}) whose constraint matrix $A$ is 
totally unimodular and all the elements of $b$ are integers.
Recall that the matrix $A$ is totally unimodular if the determinant of every nonsingular square submatrix of $A$ is 1 or -1.
Then all the elements of any BFS are integers, so $\delta\ge 1$. 
Let us bound $\gamma$. 
Let $(x_B,0)\in\Re^{m}\times \Re^{n-m}$ be a basic feasible solution of (\ref{primal}). 
Then we have $x_B=A_B^{-1}b$. 
Since $A$ is totally unimodular, all the elements of $A_B^{-1}$ are $\pm 1$ or $0$.
Thus for any $j\in B$ we have $x_j \le \|b\|_1$, which implies that $\gamma\le \|b\|_1$. 
By Theorem \ref{iteration}, we obtain the following result. 
\begin{corollary}
Assume that the constraint matrix $A$ of (\ref{primal}) is totally unimodular and the constraint vector 
$b$ is integral. 
When we apply the simplex method with the most negative rule or the best improvement rule for (\ref{primal}), 
we encounter at most $n\lceil m\|b\|_1\log (m\|b\|_1)\rceil$ different basic feasible solutions. 
\end{corollary} 
\subsection{Markov Decision Problem}
The Markov Decision Problem (MDP), where the number of possible actions is two, is  formulated as 
\begin{equation}
\begin{array}{ll}
\min &c_1^Tx_1+c_2^Tx_2,\\
{\rm subject~to}&(I-\theta P_1)x_1+(I-\theta P_2)x_2=e,\\ \label{mdp}
&x_1,x_2\ge 0,
\end{array}
\end{equation}
where $I$ is the $m\times m$ identity matrix, $P_1$ and $P_2$ are $m\times m$ Markov matrices,
$\theta$ is a discount rate, and $e$ is the vector of all ones. 
MDP(\ref{mdp}) has the following properties. 
\begin{enumerate}
\item MDP(\ref{mdp}) is nondegenerate. 
\item The minimum value of all the positive elements of BFSs is greater than or equal to 1, or 
equivalently, $\delta\ge 1$.
\item The maximum value of all the positive elements of BFSs is less than or equal to $\frac{m}{1-\theta}$, 
or equivalnetly, $\gamma \le \frac{m}{1-\theta}$.
\end{enumerate}
Therefore we can apply Corollary \ref{iteration} and obtain a similar result to Ye \cite{ye}.
\begin{corollary}The simplex method for solving MDP (\ref{mdp}) finds an optimal solution in at most 
$n\lceil \frac{m^2}{1-\theta}\log\frac{m^2}{1-\theta}\rceil $ iterations, where $n=2m$.
\end{corollary}
\vspace{5mm}
{\bf Acknowledgement}\\
This research is supported in part by Grant-in-Aid for Science Research (A) 
20241038 of Japan Society for the Promotion of Science.


\begin{thebibliography}{99}
\bibitem{dantzig} G. B. Dantzig: 
{\it Linear Programming and Extensions}. 
Princeton University Press, Princeton,
New Jersey, 1963.
\bibitem{km} V. Klee and G. J. Minty: 
How good is the simplex method. In O. Shisha, editor, 
{\it Inequalities III}, Academic Press, New York, NY, 1972.
\bibitem{ye} Y. Ye: The Simplex Method is Strongly Polynomial for the
Markov Decision Problem with a Fixed Discount Rate. 
Technical paper, available at http://www.stanford.edu/~yyye/simplexmdp1.pdf, 2010.
\end{thebibliography}
\end{document}